\documentclass{article}
\usepackage{amsmath,amsthm,amssymb}

\newcommand{\fc}{\mathfrak{c}}

\newcommand{\ga}{\alpha}
\newcommand{\gb}{\beta}

\newcommand{\gd}{\delta}
\newcommand{\gw}{\omega}
\newcommand{\eps}{\epsilon}

\newcommand{\gS}{\Sigma}
\newcommand{\gs}{\sigma}

\newcommand{\R}{\mathbb{R}}
\newcommand{\cP}{\mathcal{P}}
\newcommand{\cB}{\mathcal{B}}

\newcommand{\cof}{\mathtt{cof}}

\newcommand{\lnull}{\mathtt{null}}

\newcommand{\lr}{L(\mathbb R)}

\newcommand{\rng}{\mathrm{rng}}

\newcommand{\pioneoneonsigmaoneone}{$\mathbf{\Pi^1_1}$ on $\mathbf{\Sigma^1_1}$}

\newtheorem{theorem}{Theorem}[section]
\newtheorem{lemma}[theorem]{Lemma}
\newtheorem{claim}[theorem]{Claim}
\newtheorem{corollary}[theorem]{Corollary}

\newtheorem{proposition}[theorem]{Proposition}

\theoremstyle{definition}
\newtheorem{definition}[theorem]{Definition}
\newtheorem{example}[theorem]{Example}

\newcommand{\cantor}{2^\gw}
\newcommand{\baire}{\gw^\gw}

\newcommand{\liff}{\leftrightarrow}

\newcommand{\power}{\mathcal P}
\newcommand{\gwtree}{\gw^{<\gw}}

\title{Two preservation theorems}

\author{
Jind{\v r}ich Zapletal
\thanks{Partially supported by GA {\v C}R grant 201-03-0933 and NSF grant DMS 0300201.}\\
University of Florida}

\begin{document}
\bibliographystyle{plain}
\maketitle
\begin{abstract}
I prove preservation theorems for countable support iteration of proper forcing concerning certain classes
of capacities and submeasures. Various examples of forcing notions motivated by measure theory are included.
\end{abstract}

\section{Introduction}
I will use the ideal approach to definable proper forcing \cite{z:book} to prove two preservation theorems
for the countable support iteration. Neither seems to have a counterpart in the classical combinatorial
approach to iterated forcing, such as \cite{bartoszynski:set}. It is my hope that the comparison of the
two methods will bring deeper understanding of the subject.

The first preservation theorem deals with strongly subadditive capacities. Recall

\begin{definition}
\cite{kechris:classical}, 30.1.
Let $X$ be a Polish space. A \emph{capacity} on $X$ is a function $c:\cP(X)\to\R^+\cup\{\infty\}$ such that

\begin{enumerate}
\item $c(0)=0$, $A\subset B\to c(A)\leq c(B)$
\item $c$ is continuous in countable increasing unions: if $A_0\subset A_1\subset\dots\subset X$ are sets then
$c(\bigcup_n A_n)=\sup_nc(A_n)$
\item for a compact set $C\subset X$, $c(C)=\inf\{c(O):C\subset O$ and $O\subset X$ is an open set$\}<\infty$.
\end{enumerate}

The capacity $c$ is \emph{outer regular} if $c(A)=\inf\{c(O):A\subset O$ and $O\subset X$ is an open set$\}$
for all sets $A\subset X$. The capacity $c$ is \emph{strongly subadditive} if $c(A\cup B)+c(A\cap B)\leq
c(A)+c(B)$ holds for all sets $A,B\subset X$.
\end{definition}

Outer regular capacities are key objects in measure theory 
\cite{rogers:hausdorff} and potential theory \cite{adams:potential}. They also generate
a canonical family of proper forcings \cite{z:newton}. The basic examples
are the outer Lebesgue measure and the Newtonian capacity \cite{adams:potential}. 
Most if not all capacities in potential theory 
are strongly subadditive.

All strongly subadditive capacities are outer regular.
Outer regular capacities are determined by their values on any basis for the Polish space $X$
closed under finite unions. Such a definition of a capacity can then be used in various models of set theory.
Obvious absoluteness problems arise:

\begin{definition}
Suppose that $P$ is a forcing and $c$ is an outer regular capacity on a Polish space $X$. The forcing
$P$ \emph{preserves} the capacity $c$ if for every set $A\subset X$, $P\Vdash\dot c(\check A)=\check c(\check A)$.
\end{definition}

\begin{theorem}
\label{firsttheorem}
(LC) Suppose that $c$ is a strongly subadditive capacity. 
Suppose that $P$ is a proper universally Baire real forcing. If $P$ preserves the capacity $c$
then the countable support iterations of $P$ of any length preserve the capacity $c$ as well.
\end{theorem}

Here, LC denotes a suitable large cardinal assumption which implies determinacy of all games
used in the argument. An inaccessible
cardinal $\gd$ which is a limit of $<\gd$ strong cardinals and Woodin cardinals will suffice \cite{z:book},
\cite{neeman:book}.
By a universally Baire proper real forcing I mean a forcing satisfying Definition 2.1.5
of \cite{z:book}. This includes most proper forcings
defined from real parameters for adding a single real. Note that the statement of the theorem
abstractly implies its version in which there are more than just one iterand and the iterands 
alternate in some predictable manner.

The argument immediately extends beyond the class of strongly subadditive capacities. In fact, it
is enough for the capacity to be outer regular and to be continuous under increasing wellordered unions
under a hypothesis such as AD+. This includes Hausdorff capacities \cite{rogers:hausdorff} as well as capacities implicitly
constructed by Stepr\=ans in \cite{steprans:many}, among others. I do not know an example of an outer
regular capacity which would not be continuous in this way. The last section of the present paper
indicates a determinacy argument that can prove continuity for many capacities.

The second preservation theorem deals with Method I or pavement submeasures. Recall

\begin{definition}
\cite{munroe:measure}
A \emph{submeasure} on a Polish space $X$ is a function $\mu:\power(X)\to\R^+$ with the properties
$\mu(0)=0, A\subset B\to\mu(A)\leq\mu(B)$ and $\mu(A\cup B)\leq\mu(A)+\mu(B)$. A \emph{pavement submeasure}
is one defined from a collection $U$ of \emph{pavers}, subsets of the space $X$, and a 
\emph{weight function} $w:U\to\R^+$ extended to $w:\power(U)\to\R^+$ by $w(V)=\gS_{u\in V}w(u)$,
by setting $\mu(A)=\inf\{w(V):V\subset U$ is a set such that $A\subset\bigcup V\}$.
\end{definition}

Submeasures defined by a countable collection of Borel pavers are frequently 
encountered in measure theory \cite{rogers:hausdorff}.
They also generate canonical proper forcings \cite{z:four}. Such submeasures can be identified with their definitions,
and an obvious absoluteness problem will be encountered again. The iteration theorem now says

\begin{theorem}
\label{secondtheorem}
(LC) Suppose that $\mu$ is a pavement submeasure derived from a countable system of Borel pavers. Suppose that
$P$ is a universally Baire proper real forcing. If $P$ preserves the submeasure $\mu$ then so do all of its
countable support iterations.
\end{theorem}

Even though it is easier to define a pavement submeasure than a capacity, the pavement submeasures are not
encountered as frequently in practice. Some of them are in fact capacities, such as the outer Lebesgue measure,
or the Hausdorff capacities. Others are very far from capacities, such as the following submeasure
$\mu$ connected with $\baire$-bounding. Consider the Baire space $X=\baire$ with the countable collection
$U=\{u_{t,n}:t\in\gwtree, n\in\gw\}$ of pavers where $u_{t,n}=\{f\in\baire:t\subset f, f(|t|)\in n\}$,
and a weight function $w:U\to\R^+$ given by $w(u_{t,n})=\eps_t$ where $\eps_t:t\in\gwtree$ are positive numbers
with finite sum. It is not difficult to see that a forcing is $\baire$-bounding if and only if it preserves
the submeasure $\mu$. Thus the theorem subsumes the known fact about iteration of bounding forcings, at least
for the case of definable forcings.

I do not know if the theorems can be proved without large cardinal assumptions or without the restriction
to the class of definable forcings. Note that \cite{z:book}, Theorem 5.4.16 gives a natural preservation theorem which fails
in the context of undefinable forcings. In terms of relative consistency, the statements of the theorems are not strong.
All properties of the capacities or submeasures needed are in fact satisfied in
the Levy collapse model derived from a single inaccessible cardinal, and if one is willing to start his
iteration from that model, the preservation theorems will hold as long as the iteration is definable. 

The notation of this paper sticks to the set theoretic standard of \cite{jech:set}. 
As the standard reference for definable forcing I use \cite{z:book}, for measure theory \cite{rogers:hausdorff},
for descriptive set theory \cite{kechris:classical}, and for cardinal invariants \cite{bartoszynski:set}.
If $I$ is an ideal on
some set $X$, the statement $\forall^Ix\ \phi(x)$ is a shorthand for $\{x\in X:\lnot\phi(x)\}\in I$.
The letter $\lambda$ is reserved for the outer Lebesgue measure. A set is capacitable for a given capacity
if it has compact subsets of capacity arbitrarily close to its own.

\section{Capacities}

The key property of strongly subadditive capacities is stated in the following:

\begin{lemma}
(LC) Suppose that $c$ is a strongly subadditive capacity.
\begin{enumerate}
\item Every universally Baire set has a Borel subset of equal capacity.
\item If $A\subset\R$ is a universally Baire set then $\lr[A]\models c$ is continuous in wellordered increasing unions.
\end{enumerate}
\end{lemma}

The lemma will be proved in the last section of this paper. 
It is exactly here where it is necessary to restrict attention to the definable context. I do not know if it is
consistent with ZFC that strongly subadditive capacities are continuous in increasing unions of length $\gw_1$
consisting of Borel sets.

Towards the proof of Theorem~\ref{firsttheorem}, suppose that $P$ is a universally Baire proper real real forcing in the sense
of Definition 2.1.5 of \cite{z:book}. 
By Lemma 2.1.6 there, there is a $\gs$-ideal $I$ on some Polish space $Y$ such that $P$
is in the forcing sence equivalent to $P_I$, the poset of Borel sets modulo the ideal
$I$, every universally Baire set is either in the ideal $I$ or has an $I$-positive Borel subset, and the collection
of all Borel sets in the ideal $I$ is universally Baire.

Suppose that $c$ is an outer regular capacity on some Polish space $X$. 

\begin{lemma}
\label{helplemma}
(LC) The following are equivalent:
\begin{enumerate}
\item The forcing $P$ preserves capacity
\item for every number $\eps>0$, every $I$-positive Borel set $C\subset Y$ and every
universally Baire set $B\subset X\times C$, if every horizontal section $B^y$ has capacity $\leq\eps$ then
the set $B^I=\{x\in X:\forall^I y\in C\ \langle y,x\rangle\in B\}\subset X$ has capacity $\leq\eps$.
\end{enumerate}
\end{lemma}

\begin{proof}
If (2) failed with some $C\subset Y$, $B\subset X\times C$ and $\eps$, then the 
set $B^I\subset X$ has capacity $>\eps$, but in the extension the set $B^I\cap V$ is forced to be covered
by the generic section of the set $B$, which by the universally Baire absoluteness has capacity $\leq\eps$,
showing that the capacity of the set $B^I$ was not preserved.

On the other hand, if (1) failed then there would be a set $A\subset X$, a condition $C\in P_I$ and a $P_I$-name $\dot O$ for
an open set of capacity $<c(A)$ such that $C\Vdash\check A\subset\dot O$. By a usual uniformization argument
we may assume that there is a Borel set $B\subset X\times C$ such that every horizontal section of it
is an open set of capacity $<c(A)$ and $C$ forces $\dot O$ to be the generic section of the set $B$. Now
note that necessarily $A\subset B^I$ and so $c(B^I)\geq c(A)$ and (2) fails.
\end{proof}

I will argue for Theorem~\ref{firsttheorem} in the case of iterations of length $\gw$, the rest is
routine. I need to recall some basic facts about countable support iteration as they appear in \cite{z:book}.
The ideal $I^\gw$ on the space $Y^\gw$ consists of those sets $A\subset Y^\gw$ for which Adam has a winning strategy
in the game $G_\gw(A)$ of length $\gw$ in which he plays sets $A_n:n\in\gw$ in the ideal $I$ and Eve counters with points
$y_n\notin A_n:n\in\gw$; Eve wins if the sequence of her answers belongs to the set $A$. It turns out
that

\begin{itemize}
\item the countable support iteration of length of the poset $P_I$ is in the forcing sense isomorphic to the poset
$P_{I^\gw}$
\item every universally Baire set $A\subset Y^\gw$ is either in the ideal $I^\gw$ in which case Adam even has
a universally Baire winning strategy, or it has a Borel $I$-perfect subset where
\item a Borel set $C\subset Y^\gw$ is Borel $I$-perfect if there is a tree $p\subset Y^{<\gw}$ such that
$p\cap Y^n$ is Borel for all $n\in\gw$, every node splits into $I$-positively many immediate successors and
$B$ is the set of all cofinal branches in the tree $p$.
\end{itemize}

From now on I will identify conditions in the iteration $P^\gw$ with trees as in the last item. Note 
that every universally Baire $I$-positively branching tree contains a condition in the iteration as a subset by the
second item. For a condition $p\in P^\gw$ write $[p]\subset Y^\gw$ for the set of all cofinal branches through the tree $p$.

Suppose $p$ is a condition in the iteration, and $p\Vdash\dot A\subset\dot X$ is a set of capacity $\leq\eps$.
We must show that the set $\{x\in X:p\Vdash\check x\in\dot A\}$ has capacity $\leq\eps$. 
By the outer regularity, we may assume that $\dot A$ is forced to be an open set,
$\dot A=\bigcup\rng\dot f$ for some function $\dot f$ from $\gw$ to some fixed basis $\cB$ for the Polish space $X$. 
By a standard
countable support iteration fusion argument, there is a condition $q\leq p$ and a Borel function $F:[q]\to\cB^\gw$
such that

\begin{itemize}
\item for every sequence $\vec y\in [q]$, $c(\bigcup\rng\ F(\vec y))\leq\eps$
\item $q\Vdash\dot f=\dot F(\vec y_{gen})$ where $\vec y_{gen}$ is the name for the generic $\gw$-sequence
of points in the Polish space $Y$
\item for every number $n$ and every sequence $\vec y\in[q]$, the value $F(\vec y)(n)$ depends only on
$\vec y\restriction n$; in other words, the value $\dot f(n)$ has been decided at the $n$-th stage of the iteration.
\end{itemize}

In view of the last item I will extend the function $F$ to act on the nodes in the tree $q$ so that
$F(\vec y)\subset\cB$ is the set of the unique values $F(\vec z)(m):m\in|\vec y|$ for any infinite sequence $\vec z\in [q]$
extending the finite sequence $\vec y$. So $\rng\ F(\vec z)=\bigcup\{F(\vec y):\vec y\in q,\vec y\subset\vec z\}$ for every
sequence $\vec z\in[q]$.

It will be enough to show that the set $\{x\in X:q\Vdash\check x\in\dot A\}$ has capacity $\leq\eps$.
By induction on an ordinal $\ga$ for all sequences $\vec y\in q$
simultaneously define sets $A(\ga,\vec y)$ in the following fashion:

\begin{itemize}
\item $A(0,\vec y)=\bigcup F(\vec y)$
\item $A(\ga+1,\vec y)=\{x\in X:\forall^J z\in Y\ \vec y^\smallfrown z\in q\to x\in A(\ga,\vec y^\smallfrown z)\}$
\item $A(\ga,\vec y)=\bigcup_{\gb\in\ga}A(\gb,\vec y)$ for limit ordinals $\ga$.
\end{itemize}

By simultaneous induction on $\ga$ it is easy to prove that $\vec y\subset\vec z\to A(\ga, \vec y)\subset
A(\ga, \vec z)$ and $\gb\in\ga\to A(\gb,\vec y)\subset A(\ga,\vec y)$. 
Note that the whole process takes place in a model such as $\lr[I]$ all of whose sets of reals are universally Baire and 
which satisfies AD+.
This means that by transfinite induction again it is possible to argue that $c(A(\ga,\vec y))\leq\eps$ for
every finite sequence $\vec y\in q$ and every ordinal $\ga$.
Just use the assumption that the forcing $P_I$ preserves capacity at the successor stages
of the induction, and the continuity of the capacity in increasing wellordered unions at the limit
stages.

The inductive process must stabilize at some ordinal $\Omega<\fc^+$. I claim that
the set $\{x\in X:B\Vdash\check x\in\dot A\}$ is included in $A(\Omega,0)$, which will complete the proof
given the fact that $c(A(\Omega,0))\leq\eps$. Suppose that $x\notin A(\Omega,0)$. The tree
$r=\{\vec y\in q:x\notin A(\Omega,\vec y)\}\subset q$ is universally Baire and $I$-positively branching, 
therefore it contains a condition $s\in P^\gw$ as a subset. It is immediate from the definitions
that $s\Vdash\check x\notin\dot A$ as desired. Theorem~\ref{firsttheorem} follows.

An inquisitive mathematician will ask, how long is the transfinite process? Perhaps in some cases it is short enough
so that I can remove the determinacy argument about continuity in increasing unions? It turns out that
in general no clear answer is available. I know of only one special case in which it is possible to eliminate the
large cardinal hypothesis entirely.

\begin{theorem}
\label{thirdtheorem}
Suppose that $I$ is an iterable, \pioneoneonsigmaoneone\ ideal on some Polish space $Y$. If the forcing
$P_I$ preserves outer Lebesgue measure then so do all of its countable support iterations.
\end{theorem}

The key assumption here is the restriction to the \pioneoneonsigmaoneone\ ideals. 
This is a nontrivial restriction--the Miller forcing is \pioneoneonsigmaoneone\
while the Laver forcing is not \cite{z:book} Sections 2.3.2 and 2.3.6. 
There are definable forcings which are not \pioneoneonsigmaoneone\
and preserve many outer regular capacities; the Laver forcing is one example.
The iterability condition on an ideal is essentially the demand that the poset
$P_I$ is provably proper--see Definition 3.1.1 of \cite{z:book}.
I also need to deal with the outer Lebesgue measure in order to
be assured of measurability of coanalytic sets. Coanalytic sets are capacitable for all strongly subadditive
capacities under quite mild large cardinal assumptions, but I do not know if they are capacitable in ZFC.

\begin{proof}
The argument is the same as for Theorem~\ref{firsttheorem}, paying attention to the definability of the sets
encountered, and in the end using a theorem of Cenzer and Mauldin which puts a bound on the length of the transfinite
sequence.

First note that whenever $B\in P_I$ is an $I$-positive Borel set and $D\subset B\times [0,1]$ is a coanalytic set
all of whose vertical sections have Lebesgue measure $\leq\eps$ then the set $C=\{r\in [0,1]:\forall^I y\in B\ 
\langle y,r\rangle\in D\}$ is coanalytic of Lebesgue measure $\leq\eps$. First, $C$ is coanalytic since
the ideal $I$ is \pioneoneonsigmaoneone. And second, the set $C$ consists only of reals which are forced by $B$ to
belong to the generic vertical section of the set $D$. The fact that all vertical sections of the set $D$ have Lebesgue
measure $\leq\eps$ is $\mathbf{\Pi^1_2}$ ($\forall y\in B\ \forall C$ compact of Lebesgue measure
$>\eps\ \exists r\in C\ r\notin D_y$) and therefore absolute between the ground model and the generic extension.
Since the forcing $P_I$ preserves the outer measure of the set $C$, it must be that $\lambda(C)\leq\eps$.

Continue in the line of the argument for Theorem~\ref{firsttheorem}, using the \pioneoneonsigmaoneone\ property
of the ideal and results of \cite{z:book} Section 3.3 to set up the iteration in ZFC. In order to rein in the transfinite
process in the end of the argument, given a condition $p$ in the iteration consider the space
$Z=p\times [0,1]$ and the coanalytic inductive monotone operator $\Gamma$ on it given by $\langle s,r\rangle\in\Gamma(A)$
iff $\langle s,r\rangle\in A$ or $\forall^I y\in Y\ s^\smallfrown y\in p\to\langle s^\smallfrown y,r\rangle\in A$.
It is clear that the transfinite process indicated in the proof of Theorem~\ref{firsttheorem}
is just the iteration of the operator $\Gamma$
on a certain Borel initial set $A$. A theorem of Cenzer and Mauldin \cite{cm:inductive} 1.6(e)
shows that the iteration stabilizes after $\leq\gw_1$
many steps in a coanalytic set $A_{\gw_1}$, and every analytic subset of $A_{\gw_1}$ is actually a subset
of some countable iterand. Consider the set $C=\{r\in [0,1]: \langle 0,r\rangle\in A_{\gw_1}\}$. As before,
the set $C$ includes every real $r\in [0,1]$ such that $p\Vdash\check r\in\dot O$: if $r\notin C$ then the tree
$\{s\in p: \langle s,r\rangle\notin A_{\gw_1}\}\subset p$ is analytic and $I$-positively branching, by
Lemma 3.3.1 of \cite{z:book} it contains a condition in the iteration as a subset, which then forces
$\check r\notin\dot O$. Also, $\lambda(C)\leq\eps$: if this failed, the set $C$ would have a compact subset
of measure $>\eps$; however, such a compact set has to be exhausted at some countable stage of the iteration,
an impossibility since the sets appearing at countable stages of the transfinite process 
all have Lebesgue measure $\leq\eps$.
\end{proof}

\section{Pavement submeasures}

Suppose that $\mu$ is a Method I submeasure on a Polish space $X$, built from a countable set $U$ of Borel pavers
and a weight function $w:U\to\R^+$. I will again need several determinacy-related basic properties of the submeasure
$\mu$.

\begin{lemma}
\label{helplemma2}
(LC) Let $\mu$ be a pavement submeasure built from a countable set of Borel pavers.
\begin{enumerate}
\item every universally Baire set has a Borel subset of the same $\mu$-submeasure
\item Player Nonempty has a winning strategy in the descending chain game with $\mu$-positive universally Baire sets.
\end{enumerate}
\end{lemma}

Here, the descending chain game is played between players Empty and Nonempty. In the beginning, 
Empty indicates a $\mu$-positive universally Baire set $B_0$ and then the players alternate to obtain
a descending chain $B_0\supset B_1\supset B_2\supset\dots$ of universally Baire $\mu$-positive sets.
Player Nonempty wins if the intersection of the chain is nonempty. This is closely related to
the precipitousness games played with other ideals, in which sets of arbitrary complexity are allowed.
The restriction to universally Baire sets here gives player Nonempty his winning strategy, and this is exactly
the reason why I do not know how to prove the theorem in undefinable context. I do not know if it is consistent
for the related null ideals to be precipitous.

Now (1) above will be proved in the last section of this paper. (2) is an immediate consequence of (1). Let
$J$ be the ideal of $\mu$-null sets and recall that the factor poset $P_J$ of Borel $J$-positive sets
adds a real which belongs to all sets in the generic filter \cite{z:book} Lemma 2.1.1. In the run of the play of the descending
chain game let the player Nonempty play just Borel sets--this is possible by the first item--and on the side
generate an increasing sequence $M_0\subset M_1\subset\dots$ of countable elementary submodels of some large structure,
and use some simple bookkeeping to ensure that the descending chain of Borel sets generates an $M$-generic filter $g
\subset P_J$,
where $M=\bigcup_nM_n$. Then $M[g]\models$the generic real belongs to all sets in the filter $g$, so the
generic real shows that player Nonempty won. A noteworthy aside--the forcing $P_J$ is in fact proper by the results
of \cite{z:four} Section 6.

Towards the proof of Theorem~\ref{secondtheorem}, suppose that $P$ is a universally Baire proper real forcing.
Repeat the setup from the previous section. There is a $\gs$-ideal $I$
on some Polish space $Y$ such that $P=P_I$, $I$ satisfies the weak dichotomy: every universally Baire set is either
in the ideal or else it has a Borel $I$-positive subset, and the collection of all Borel $I$-positive sets
is in a suitable sense universally Baire. Lemma~\ref{helplemma} deserves a restatement.

\begin{lemma}
(LC) The following are equivalent:

\begin{enumerate}
\item $P_I$ preserves the submeasure $\mu$
\item for all sets $B,C,D$ such that $B\subset X$ is universally Baire and $\mu$-positive, $C\subset Y$ is universally
Baire $I$-positive set, and $D\subset B\times C$ is a universally Baire set with $I$-small vertical sections, there
is a horizontal section of the set $(B\times C)\setminus D$ which has $\mu$ submeasure equal to the submeasure
of the whole set $B$.
\end{enumerate}
\end{lemma}

I will prove the theorem for iterations of length $\gw$, the argument for larger ordinals is then routine.
Let $\vec y_{gen}$ denote the generic $\gw$-sequence of points in the Polish space $Y$.
Suppose $p\in P_\gw$ is a condition and $\dot W$ is a $P$-name for a subset of $U$ of weight $<\eps$. 
I must show that the set $B=\{x\in X:p\Vdash\check x\in\bigcup\dot W\}$ has $\mu$-submeasure $<\eps$.
By a standard countable support iteration fusion argument just as in the previous section,
strengthening the condition $p$ if necessary I may assume that there is a Borel function
$F:p\to [U]^{\aleph_0}$ such that $\vec y\subset \vec z\in p$ implies $F(\vec y)\subset F(\vec z)$,
$w(F(\vec y))<\eps$, and $p\Vdash\dot W=\bigcup_n\dot F(\vec y_{gen}\restriction n)$. I will extend the
function $F$ to act on all infinite sequences in the set $[p]$ by $F(\vec z)=\bigcup_n F(\vec z\restriction n)$.
The set $B$ is then universally Baire: $B=\{x\in X:
A_x\in I^\gw\}$ where $A_x=\{\vec y\in [q]:x\notin\bigcup F(\vec y)\}$, 
and it is possible to find a universally Baire assignment
$x\in B\mapsto\tau_x$ where $\tau_x$ is a winning strategy for Eve in the game $G_\gw(A_x)$.

Now suppose for contradiction that $\mu(B)>\eps$.
Fix a winning strategy $\gs$ for the Nonempty player in the decreasing chain game with universally Baire
$\mu$-positive sets. Build by induction on $n\in\gw$ sequences $\vec y_n\in q$ of length $n$ and
a decreasing chain of $\mu$-positive sets $B\supset B_0\supset B_1\supset\dots$ so that

\begin{itemize}
\item $\vec y_n$ is a legal sequence of Eve's answers against the strategy $\tau_x$ for every point $x\in B_n$
\item $B_n\cap \bigcup F(\vec y_n)=0.$
\end{itemize}

Suppose this is possible to arrange so that in the end the intersection $\bigcap_n B_n$ is nonempty.
Writing $x$ for any point in the intersection and $\vec y=\bigcup_n y_n$, it should be the
case that both $\vec y\in A_x$ by the second item, and $\vec y\notin A_x$ by the first item, and this
will be the contradiction which will complete the proof of the theorem.

However, in order to perform the induction successfully it is necessary to add some extra statements to the induction hypothesis.
A bit of notation: if $\vec y\in p$ is a finite sequence, a $\vec y$-tree is a Borel tree $r\subset p$ with
trunk $\vec y$ which branches into $I$-positively many immediate successors at every node past the sequence $\vec y$.
In the course of the induction I will also build $\vec y_n$-trees $p_n\subset q$, $\mu$-positive sets $\bar B_n\subset X$
and finite sets $W_n:n\in\gw$ so that

\begin{itemize}
\item $p\supset p_0\supset p_1\supset p_2\dots$
\item $F(\vec y_n)\subset W_n$, $B_n\cap\bigcup W_n=0$ and for every sequence $\vec z\in[p_n]$ 
it is the case that $W_n\subset F(\vec z)$, the weight of
the set $F(\vec z)\setminus W_n$ is smaller than $\mu(B_n)$, and in fact the numbers 
$\{w(F(\vec z)\setminus W_n):z\in [p_n]\}$
are bounded below $\mu(B_n)$
\item $B_0\supset \bar B_0\supset B_1\supset\bar B_1\supset\dots$ and $\bar B_0, B_1,\bar B_1, B_2,\dots$
constitutes a run of the decreasing chain game respecting the strategy $\gs$.
\end{itemize}

To perform the induction, at $n=0$ let $p_0=p, W_0=F(0), B_0=B\setminus\bigcup
F(0), \vec y_0=0$. Suppose that the sequence $\vec y_n$,
tree $p_n$ and sets $W_n$ and $B_n$ have been constructed. Let $C_n=\{y\in Y:
\vec y_n^\smallfrown y\in p_n\}\notin I$. Considering the rectangle $B_n\times C_n$ and the set $D_n=\{\langle x,y\rangle
\in B_n\times C_n:y\in\tau_x(\vec y)\}\subset B_n\times C_n$, using the assumption that the poset $P$ preserves the submeasure
and using the Lemma~\ref{helplemma2}, 
conclude that there is a point $y\in C_n$ such that the set $\tilde B_n=\{x\in B_n:
y\notin\tau_x(\vec y)\}$ has measure equal to that of $B_n$. Put $\vec y_{n+1}=\vec y_n^\smallfrown y$.
To find the remaining objects, construct a decreasing sequence $q_n^m:m\in\gw$ of $\vec y_{n+1}$-trees 
below the tree $p_n\restriction\vec y_{n+1}$ such that
for each $m\in\gw$ there is a finite set $V_n^m\subset U$ such that
for every sequence $\vec z\in [q_n^m]$ it is the case that $V_n^m\subset F(\vec z)\setminus W_n$ and
the weight of the set $F(\vec z)\setminus (V_n^m\cup W_n)$ is smaller than $2^{-m}$. This is possible 
by the countable completeness of the ideal $I^\gw$. By the second item of the induction hypothesis the
weights of the sets $V^m_n\subset U$ are bounded below $\mu(\tilde B_n)$ and so $w(\bigcup_m V_n^m)<\mu(\tilde B_n)$
and the set $\bar B_n=\tilde B_n\setminus\bigcup_m V_n^m$ is $\mu$-positive. Let $B_{n+1}$ be the answer
of the strategy $\gs$ to the set $\bar B_n$, let $m$ be so large that $2^{-m}<\mu(B_{n+1})$ and let
$p_{n+1}=q_n^m$. The induction hypotheses continue to hold.

In the end, the intersection $\bigcap_n B_n$ is nonempty as desired since the strategy $\gs$ is winning for Player
Nonempty. Theorem~\ref{secondtheorem} follows.

\section{Examples}

The real reasons for which definable forcings preserve capacities such as outer Lebesgue measure are generally
neglected by forcing specialists. I decided to add a section which describes several arguments and connections.

\begin{proposition}
Suppose that $c$ is an outer regular capacity on a Polish space $X$
such that every coanalytic set is capacitable. Then Laver
forcing preserves $c$.
\end{proposition}

\begin{proof}
Similar to Theorem~\ref{thirdtheorem}. 
Suppose $T\subset\gwtree$ is a Laver tree forcing $\dot O$ to be a set of capacity $\leq\eps$.
I have to show that the set $\{x\in X:T\Vdash\check x\in\dot O\}$ has capacity $\leq\eps$.

Enumerate some basis $\cB$ of the space $X$ by $\langle V_n:n\in\gw\rangle$. Use a standard fusion argument
to find a tree $S\subset T$ such that for every node $s\in S$ the condition $S\restriction s$ decides
the statements $\dot V_n\subset\dot O$ for all $n\in |s|$. To simplify the notation assume that
$S=\gwtree$.

Consider the space $Y=\gwtree\times X$ and the operator $\Gamma:\power(Y)\to\power(Y)$ on it defined by
$\langle s,x\rangle\in\Gamma(B)\liff \langle s,x\rangle\in B\lor\forall^\infty n\ \langle s^\smallfrown n, x\rangle\in B$.
This is a monotone inductive coanalytic operator, and therefore by a theorem of Cenzer and Mauldin
\cite{cm:inductive} 1.6, given
a coanalytic set $A\subset Y$, the transfinite sequence 
given by the description $A=A_0$, $A_{\ga+1}=\Gamma(A_\ga)$ and $A_\ga=\bigcup_{\gb\in\ga}
A_\gb$ for limit ordinals $\ga$,
stabilizes at $\gw_1$ in a coanalytic set $A_{\gw_1}$ such that for every analytic set $C\subset A_{\gw_1}$
there is an ordinal $\ga\in\gw_1$ such that $C\subset A_\ga$.

Now consider the set $A\subset Y$ given by $\langle s,x\rangle \in A$ iff for some number $n\in |s|$,
$S\restriction s\Vdash\dot V_n\subset\dot O$ and $x\in V_n$. It is not difficult to see that
writing $A^s$ for the set $\{x\in X:\langle s,x\rangle\in A\}$ it is the case that $s\subset t\to A^s\subset A^t$,
these sets have capacity $\leq\eps$ and this feature persists through the countable stages of the iteration.
To see that $c(A^s_{\ga+1})\leq\eps$ note that the set $A^s_{\ga+1}$ is an increasing union of the sets
$\bigcap_{m>n}A^{s^\smallfrown m}_\ga:n\in\gw$, each of them of capacity $\leq\eps$, and use the continuity
of the capacity under increasing unions. At limit stages, use the continuity of the capacity again to argue
that $c(A^s_\ga)\leq\eps$.

Consider the coanalytic set $A_{\gw_1}$, the fixed point of the operator $\Gamma$,
and its first coordinate $B=A^0_{\gw_1}$. First note that $x\notin B$
means that $S\not\Vdash\check x\in\dot O$, since if $x\notin B$ then 
the tree $U=\{s\in S:x\notin A^s_{\gw_1}\}$ is a Laver tree
by the definition of the operator $\Gamma$ and it forces $\check x\notin\dot O$. In fact a transfinite induction argument
will show that $B=\{x\in X:S\Vdash\check x\in\dot O\}$. So it is enough to show that $c(B)\leq\eps$.
But if $c(B)>\eps$, then by the capacitability of the set $B$ there is a compact set $C\subset B$ such that
$c(C)>\eps$, and such a set must be included in the set $A^0_{\gw_1}$ for some countable ordinal
$\ga$. However, $c(A^0_{\gw_1})\leq\eps$ as proved in the previous paragraph, a contradiction!
\end{proof}

It is instructive to compare this argument with the original Woodin's proof for preservation of
outer Lebesgue measure in the Laver extension in \cite{bartoszynski:set}, 7.3.36.

It is clear that the key element in the argument was the Fubini property connecting the capacity
$c$ and the Fr\'echet ideal $J$ on $\gw$: if $\eps>0$ is a real number and $f:a\to\cP(X)$ is a function such that $a\notin J$
and $\forall n\in a\ c(f(n))\leq\eps$ then $c(\{x\in X:\forall^J n\in a\ x\in f(n)\})\leq\eps$.
It seems to be very tricky to verify the status of the Fubini property between various ideals
on $\gw$ or other sets and various capacities or submeasures.

The following general theorems concern the preservation of outer Lebesgue measure by certain classes of forcings.
In the definable context this implies that these forcings also preserve all strongly subadditive capacities,
as shown in the next proposition. I do not know an example of a forcing which preserves a given strongly subadditive
capacity such as the Newtonian capacity, and makes the set of ground model reals Lebesgue null. I have not
looked for such an example with any sincerity. Shelah and Stepr\=ans 
\cite{steprans:hausdorff}, \cite{steprans:compare} extracted a forcing which preserves a certain Hausdorff capacity
but makes the set of ground models reals Lebesgue null.

\begin{proposition}
(LC) Suppose that $P$ is a proper universally Baire real forcing. If $P$ preserves outer Lebesgue measure then it preserves
every strongly subadditive capacity.
\end{proposition}

\begin{proof}
The key ingredient is the theorem of Choquet stating that a strongly subadditive capacity is an envelope of measures:
if $c$ is a strongly subadditive capacity on a Polish space $X$ and $K\subset X$ is a compact set then
there is a Borel measure $\mu$ such that $\mu\leq c$ and $\mu(K)=c(K)$.

Fix a suitable $\gs$-ideal $I$ on a Polish space $Y$ so that $P=P_I$, every universally Baire $I$-positive set has a Borel
$I$-positive subset, and the collection of $I$-positive Borel sets is universally Baire.
Suppose that $B\in P_I$ is a condition forcing $\dot O\subset X$ is an open set of capacity $\leq\eps$.
Strengthening the condition $B$ if necessary I may assume that there is a Borel set $D\subset B\times X$
all of whose vertical sections have capacity $\leq\eps$ and $B$ forces the set $\dot O$ to be the generic vertical
section of the set $\dot O$. It will be enough to show that the set $C_0=\{x\in X:\forall^Iy\in B\ \langle y,x\rangle\in D\}
=\{x\in X:B\Vdash\check x\in\dot O\}$
has capacity $\leq\eps$.

Suppose $c(C_0)>\eps$. The set $C_0$ is universally Baire, therefore capacitable for the strongly subadditive capacity $c$
by Corollary~\ref{capacitycorollary3},
and it has a compact subset $C_1$ of capacity $>\eps$. Find a Borel measure $\mu$ on $X$
such that $\mu\leq c$ and $\mu(C_1)=c(C_1)$. Since all vertical sections of the set $D$ have $\mu$-measure $\leq\eps$,
the set $\dot O\subset X$ is forced to have $\mu$-measure $\leq\eps$. Thus
$B\Vdash\dot\mu^*(\check C_1)\leq\eps<\check \mu(\check C_1)$ and therefore the forcing $P$ does not preserve
the outer measure $\mu^*$. By the measure isomorphism theorem \cite{kechris:classical} 17.41 the forcing
$P$ does not preserve Lebesgue outer measure. Contradiction!
\end{proof}

It is possible that this Proposition is in fact true without the large cardinal assumptions and definability
restrictions if the capacitability argument can be somehow avoided. I do not know how to do that.

\begin{proposition}
Suppose that $c$ is a strongly subadditive capacity on a Polish space $X$ such that writing $I$ for the
$\gs$-ideal of zero capacity sets, the forcing $P_I$ is proper. Then $P_I$ preserves outer Lebesgue measure.
\end{proposition}

The forcing $P_I$ is proper for many strongly subadditive capacities, such as the Newtonian
capacity, as shown in \cite{z:newton}. In fact it is open whether the forcing $P_I$ is proper for
every outer regular subadditive
capacity. Forcings of this kind have not been thoroughly investigated. Since they are capacitable they
are bounding and they make the ground model reals meager
\cite{z:four} 2.17 and 7.13. The forcing associated with the Newtonian
capacity is nowhere c.c.c. A test question would be to verify whether the Newtonian forcing adds a splitting real. 

\begin{proof}
The key fact entering the argument is again
the theorem of Choquet asserting that strongly subadditive capacities are envelopes
of measures: if $c$ is strongly subadditive and $K\subset X$ is a compact set, then there is a measure
$\mu$ on $X$ such that $\mu\leq c$ and $\mu(K)=c(K)$.

Suppose that the poset $P_I$ is proper, and suppose that $B\in P_I$ is a condition forcing that $\dot O\subset [0,1]$
is an open set of outer Lebesgue measure $\leq\eps$. I must prove that the set $\{r\in [0,1]:B\Vdash\check r\in\dot O\}$
has capacity $\leq\eps$. Since the poset $P_I$ is bounding, the continuous reading of names can be applied to find
a compact $c$-positive subset $B_1\subset B_0$ and a relatively open set $D\subset B_0\times [0,1]$ 
all of whose vertical sections have outer measure $\leq\eps$ such that
$B_1\Vdash\dot O$ is the generic vertical section of the set $\dot D$. Let $C=\{r\in [0,1]:B_1\Vdash
\check r\in\dot O\}=\{r\in [0,1]:$ the $r$-th horizontal section of the complement of the set $D$ has
capacity zero$\}$. Since the set $B_1$ is compact and the set $D$ is relatively open, it is not difficult to verify
that the set $C$ is $G_\gd$. It will be enough to prove that $\lambda(C)\leq\eps$.

Let $\mu$ be a measure on the space $X$ such that $\mu\leq c$ and $\mu(B_1)=c(B_1)$. If $\lambda(C)>\eps$
then the Fubini theorem applied to the product measure $\mu\times\lambda$ says that
$\mu(B_1)\cdot\eps$ (the lower bound on the product measure of the rectangle $B_1\times C$) is less than
$\mu(B_1)\cdot\eps$ (the upper bound on the product measure of the set $D$) plus $0\cdot\eps$ (the product
measure of the complement $(B_1\times C)\setminus D$). This is a contradiction.
\end{proof}

\begin{theorem}
\label{fourththeorem}
(LC) Suppose that $\mu$ is a universally Baire
non-$\gs$-finite Borel measure on some Polish space $X$ and let $I$ be the $\gs$-ideal
generated by Borel sets of finite $\mu$-measure. Suppose that

\begin{enumerate}
\item every Borel $\mu$-positive set has a Borel subset of $\mu$-positive finite measure
\item the forcing $P_I$ is proper.
\end{enumerate}

Then the forcing $P_I$ is bounding and preserves outer Lebesgue measure.
\end{theorem}

Here a measure is universally Baire if the set $\{ \langle C,\eps\rangle: C\subset X$ is compact and $\eps=\mu(C)\}$
is universally Baire.
The point I want to make is that the investigation of property (1) above for various measures $\mu$ made
many measure theorists busy for decades--\cite{howroyd:positive}, \cite{joycepreiss:packing},
\cite{joyce:packing}, and \cite{rogers:hausdorff} Section 2.7. 

The large cardinal assumptions
can be eliminated in the case that the ideal $I$ is \pioneoneonsigmaoneone, because then the set
$C_0$ in the proof will be coanalytic and therefore Lebesgue measurable. The ideal $I$ is \pioneoneonsigmaoneone\
in all natural cases I considered. Also, the forcing $P_I$ is proper in all natural cases I considered.

\begin{example}
Let $\mu$ be the $h$-dimensional \emph{Hausdorff measure} for some compact metric space $\langle X,d\rangle$
and a continuous nondecreasing function $h:[0,\infty)\to[0,\infty)$ with $h(0)=0$. The
forcing $P_I$ is proper and bounding; for many Hausdorff measures this was proved in \cite{z:four}.
 If the function $h$ satisfies the doubling condition then property
(1) holds as proved by \cite{howroyd:positive} and the forcing $P_I$ preserves outer Lebesgue measure. In this case
the results of \cite{z:four} can be employed to show that $P_I$ does not add a splitting real. There
is a classical example of a compact metric space with a Hausdorff measure on it with no sets of finite positive measure
\cite{rogers:finite}.
I have not investigated the properties of the associated forcing.
\end{example}

\begin{example}
Let $\mu$ be the $h$-dimensional \emph{packing measure} on some compact metric space $\langle X,d\rangle$. The forcing $P_I$
is proper and bounding \cite{z:four} Section 3.3.
The status of property (1)
above depends on the precise definition of the packing measure as shown by Helen Joyce and David Preiss.
Let $h:[0,\infty)\to [0,\infty)$ be a nondecreasing continuous function with $h(0)=0$. For a set
$A\subset X$ and a real number $\gd>0$ there are three possible ways to define a $\gd$-packing
of the set $A$ and its weight:

\begin{enumerate}
\item A $1$-$\gd$-packing is a finite set of pairs $\langle x_i,r_i\rangle:i\in k$ such that
$x_i\in A$, $r_i<\gd$ and $i\neq j\to d(x_i,x_j)<r_i,r_j$; the weight of the packing is just $\gS_{i\in k}h(r_i)$
\item A $2$-$\gd$-packing is a finite set of pairs $\langle x_i,r_i\rangle:i\in k$ such that
$x_i\in A$, $r_i<\gd$ and the $r_i$-balls around $x_i$ are disjoint; the weight of the packing is just $\gS_{i\in k}h(r_i)$
\item A $3$-$\gd$-packing is a finite set of open balls with centers in the set $A$ and diameters $<\gd$;
the weight of the packing is the sum of the values of the function $h$ applied to the diameters of the balls.
\end{enumerate}

In all the three cases let the $h$-dimensional packing premeasure $\mu_p$ of the set $A$ be the infimum as $\gd>0$ of weights of
$\gd$-packings of the set $A$, and the packing measure of $A$ is the infimum of 
all numbers of the form $\gS_n\mu_p(A_n)$ where $A\subset\bigcup_n
A_n$. All the three definitions give metric measures, and they coincide in the case of ultrametric spaces.
Joyce and Preiss \cite{joycepreiss:packing} proved that the packing measures obtained through the
first definition have the property (1), and the packing measures obtained through the second
definition have the property (1) as long as the function $h$ satisfies the doubling condition.
Later Joyce \cite{joyce:packing} showed that for every function $h$ there is a compact metric
space without subsets of finite positive $h$-dimensional diameter-based packing measure.

To construct a thematic forcing which increases the invariant $\cof(\lnull)$ and keeps all other
Cicho\'n invariants unchanged (in particular preserves outer Lebesgue measure) choose a decreasing sequence of numbers
$0<d_n<1:n\in\gw$ converging to zero such that $d_n>n^2d_{n+1}$, and let $k_n:n\in\gw$ be positive natural
numbers such that the sequence $k_n\cdot d_n:n\in\gw$ diverges to infinity. Consider the metric space $X=\Pi_nk_n$
with the least difference metric $d(x,y)=d_{\Delta(x,y)}$, and consider the $1$-dimensional packing measure
$\mu$ associated with it. The numbers $k_n$ were chosen so that the space $X$ is not a countable union
of sets of finite $\mu$-measure, and the numbers $d_n$ were chosen so that any tunnel $T=\Pi_n a_n$
where $a_n\subset k_n$ is a set of size $\leq n^2$, has finite packing measure. Let $I$ be the
$\gs$-ideal on the space $X$ generated by the sets of finite packing measure. It adds an element of the space
$X$ which cannot be enclosed by an $n^2$ tunnel from the ground model, and therefore it increases
the invariant $\cof(\lnull)$ by \cite{bartoszynski:set} 2.3.9. Outer Lebesgue measure 
is preserved by the Theorem, and the forcing is bounding
and preserves nonmeager sets by the results of \cite{z:four} Section 3.3..
\end{example}

\begin{example}
Let $X,d,\nu$ be a compact metric space with 
a finite Borel measure $\nu$, and let $h:[0,\infty)\to [0,\infty)$ be a continuous nondecreasing
function such that $h(0)=0$. The $h$-dimensional \emph{Minkowski content} $\mu_p(A)$ of a set $A\subset X$
is defined as $\limsup_{\eps\to 0}\nu(A_\eps)/h(\eps)$ where $A_\eps$ is the $\eps$-neighborhood of the set $A$.
It is possible to associate a metric \emph{Minkowski measure} $\mu$ with it by the formula $\mu(A)=\inf\{\gS_n\mu_p B_n:
A\subset\bigcup_n B_n\}$. Let $I$ be the $\gs$-ideal generated by the sets of finite Minkowski measure.
Using the fact that
the closure of a set has the same Minkowski content as the set itself, it is not difficult to see
that the ideal $I$ is $\gs$-generated by a $\gs$-compact family of compact sets
and so the forcing $P_I$ is proper and bounding, \cite{z:four} Section 3.3. I do not know if every set of positive Minkowski
measure must have a subset of finite positive Minkowski measure.
\end{example}

\begin{proof}
Towards the proof of Lebesgue measure preservation part of Theorem~\ref{fourththeorem}, 
suppose that $B\in P_I$ is a condition forcing that $\dot O\subset [0,1]$ is an open set of Lebesgue measure $\leq\eps$.
By a standard properness argument, strengthening the condition $B$ if necessary I may assume that there is a Borel set
$D\subset B\times [0,1]$ such that all vertical sections of the set $D$ are open sets of Lebesgue measure
$\leq\eps$ and $B\Vdash\dot O=$the generic section of the set $\dot D$. Let $C_0=\{r\in [0,1]:\forall^I x\in B\ 
\langle x,r\rangle\in D\}$. The set $C_0$ is universally Baire and Lebesgue measurable by our large cardinal assumption,
and $C_0=\{r\in [0,1]:B\Vdash\check r\in\dot O\}$. It will be enough to show that $\lambda(C_0)\leq\eps$. Suppose
for contradiction this is not true.

The first step in the argument is to massage the sets $B_0, C_0, D$ to find an $I$-positive Borel set
$B_1\subset B_0$, a compact Lebesgue-positive set $C_2\subset C_0$, and a Borel set
$D_n\subset B_1\times C_2$ such that its vertical sections have Lebesgue measure $\leq\eps'$
for some fixed $\eps'<\lambda(C_2)$, and the horizontal section of its complement have finite $\mu$-measure
less than some fixed number $n\in\gw$. This is not difficult to do.
Choose a compact set $C_1\subset C_0$ with $\lambda(C_1)>\eps$. Note that $C_1$ forces in the random
forcing that the generic horizontal section of the complement of the set $D\subset B\times C_0$ is
$\mu$-$\gs$-finite. A standard uniformization argument with random forcing will then yield
a compact
set $C_2\subset C_1$ with $\lambda(C_2)>\eps$ and Borel sets $D_n\subset B\times C_2:n\in\gw$ such that
$D_0\supset D_1\supset\dots$, $D\cap B\times C_2=\bigcap_n D_n$ and the horizontal sections of the complement of the set $D_n$
have $\mu$-measure $\leq n$. Now for every point $x\in B_0$, the vertical sections $(D_n)_x:n\in\gw$ form a decreasing
sequence whose intersection has Lebesgue measure $\leq\lambda(D_x)\leq\eps$, and therefore there is
a number $n\in\gw$ and a condition $B_1\subset B_0$ in the poset $P_I$ such that all vertical
sections $(D_n)_x:x\in B_1$ have Lebesgue measure $\leq\eps'$ for some fixed $\eps'<\lambda(C_2)$.

Now let $m\in\gw$ be so large that $m(\lambda(C_2)-\eps')>n\lambda(C_2)$. Use property (1)
from the assumptions to find a Borel set $B_2\subset B_1$ of finite $\mu$-measure greater than $m$.
The Fubini theorem applied to the rectangle $B_2\times C_2$ and the Borel set $D_n\cap B_2\times C_2$
says that $\mu(B_2)\lambda(C_2)$ (the product measure of the rectangle) is less or equal to
$\mu(B_2)\eps'$ (an upper bound on the product measure of the set $D_n\cap (B_2\times C_2)$) plus
$n\lambda(C_2)$ (an upper bound on the product measure of the set $(B_2\times C_2)\setminus D_n$).
But this contradicts the choice of the number $m$.

To show that the forcing $P_I$ is bounding, suppose for contradiction this fails. Consider the Laver ideal
$J$ on the Baire space $\baire$, generated by all sets $A_g=\{f\in\baire:f(n)\in g(f\restriction n)\}$
as $g$ varies over all functions from $\gwtree$ to $\gw$. The Laver forcing naturally densely embeds into
the poset $P_J$--every analytic $J$-positive set contains all branches of some Laver tree, and under AD this extends
to all $J$-positive subsets of the Baire space \cite{z:book} Section 2.3.6. 
The fact that $P_I$ is not bounding is equivalent to
$I\perp J$: there is a Borel $I$-positive set $B\subset X$ and a Borel set $D\subset B\times\baire$
such that the vertical sections of the set $D$ are $J$-small and the horizontal sections of the complement
of $D$ are $I$-small \cite{z:four} 2.14. 
As in the second paragraph of the present proof, a standard uniformization argument
with the Laver forcing yields a Borel $J$-positive set $C_0\subset\baire$ and Borel sets $D_n\subset (B\times C_0):n\in\gw$
such that the horizontal sections of the complement of the set $D_n$ have $\mu$-measure $\leq n$ and
$D\cap (B\times C_0)=\bigcap_n D_n$. The rest of the argument for the preservation of Lebesgue measure
apparently has no counterpart in this situation and it is necessary to reconsider.

There are two distinct cases--either the forcing $P_I$ below the condition $B$ satisfies c.c.c. or not.
The c.c.c. case is easy to handle. Essentially by a theorem of Shelah \cite{bartoszynski:set} 3.6.47, if
the c.c.c. forcing $P_I$ below $B$ adds an unbounded real, then it adds a Cohen real, and so it makes the
the set of the ground model reals Lebesgue null. However, this is impossible by the first part of the present proof.
The non-c.c.c. case is harder. If there is an uncountable antichain in the poset $P_I$ below the condition $B$,
there must be an uncountable antichain consisting of mutually disjoint Borel sets. Each of these sets has
a Borel subset of finite positive $\mu$-measure by the property (1). Thus there is a collection
$\{ B_\ga:\ga\in\gw_1\}$ of mutually disjoint Borel subsets of the set $B$, each of positive finite $\mu$-measure.

In the Laver forcing extension, let $E_n\subset X$ be the generic horizontal section of the complement of the set $D_n
\subset (B\times C_0)$. An absoluteness argument shows that $\mu(E_n)\leq n$ and therefore the set $E_n$
has $\mu$-null intersection with all but countably many sets $B_\ga:\ga\in\gw_1$. So there is a countable
ordinal $\ga\in\gw_1$ such that for every larger countable ordinal $\gb$ it is the case that $\mu(B_\gb\cap
\bigcup_n E_n)=0$. Back to the ground model. Choose a countable elementary submodel $M$ of a large enough
structure containing all the instrumental objects, and let $C_1\subset C_0$ be the set of all $M$-generic Laver
reals. Since Laver forcing is proper, this is a $J$-positive Borel set, and the previous argument shows
that for every element $r\in C_1$ there is an ordinal $\ga_r\in\gw_1\cap M$ such that for every larger
countable ordinal $\gb\in M$ the set $B_\gb$ has $\mu$-null intersection with the $r$-th section of the complement
of the set $D$. Since there are only countably many ordinals in the model $M$, for a $J$-positive Borel
set $C_2\subset C_1$ the ordinal $\ga_r$ is the same for all elements $r\in C_2$. Choose a larger countable
ordinal $\gb\in M$, and look at the set $D\cap (B_\gb\times C_2)$ as a subset of the rectangle $B_\gb\times C_2$.
It has $J$-small vertical sections, and the horizontal sections of its complement are $\mu$-null. Thus
the condition $C_2$ in the Laver forcing forces the set $\check B_\gb$ to be covered by the generic section
of the complement of the set $D$, and therefore to have zero $\mu$ measure. However, Laver forcing preserves
outer measure, contradiction!
\end{proof}

Note how the statement and the proof of the previous Theorem tiptoe around two very unlikely issues.
I cannot exclude the possibility that the forcing $P_I$ is c.c.c.--in that case it would have to be the
Solovay forcing by the results of \cite{z:maharam}. I also cannot exclude the possibility that the
forcing $P_I$ is not c.c.c. and still has no perfect antichain consisting of mutually disjoint sets. Such
an antichain exists in all nowhere c.c.c. forcings in which an analytic family of closed sets is dense, but
I do not know if this must be the case for the forcing $P_I$. It still seems to be an open problem whether
the $\gs$-ideal of $\gs$-finite sets for a Hausdorff measure can be c.c.c., but this would have to happen
in a situation where the property (1) fails--\cite{rogers:hausdorff} Theorem 59.

\begin{example}
Let $Y$ be a Polish space and $f:Y\to\cantor$ be a Borel function which cannot be decomposed into
countably many continuous functions. Let $I$ be the ideal $\gs$-generated by the sets $A\subset Y$
such that $f\restriction A$ is continuous. The forcing $P_I$ is proper. It turns out that
the forcing $P_I$ preserves outer Lebesgue measure. This is proved by an argument essentially identical
to Example 5.4.11 of \cite{z:book}.
\end{example}

\begin{example}
Let $X,d$ be a compact metric space and let $I$ denote the ideal of $\gs$-porous subsets of the space $X$
\cite{zz:porous}.
Here, the porosity of a set $A\subset X$ at a point $x\in X$ is
$\limsup_{\gd\to 0}\sup\{r/\gd:$there is a ball of radius $r$ inside the ball around the point $x$ of radius $\gd$
which is disjoint from the set $A\}$, a set is porous if it has nonzero porosity at all of its points,
and a set is $\gs$-porous if it can be decomposed into countably many porous sets.
The forcing $P_I$ is proper by the results of \cite{z:book} Section 2.3.12, and it is bounding by the results of \cite{zz:porous}.
There is an independent determinacy argument that compact sets are dense in the poset $P_I$ 
at least in the case of zero-dimensional compact metric spaces due to Diego Rojas.
I do not know if the forcing $P_I$ in general must make the set of the ground model reals null, but I can construct an example
in which this happens. 

Let $0<p_n<1:n\in\gw$ be real numbers such that $\Pi_n p_n\neq 0$, and let $k_n:n\in\gw$
be positive natural numbers such that $\Pi_n p_n^{k_n}=0$. Consider the space $X=\Pi_n k_n$
with the least difference metric $d(x,y)=2^{-\Delta(x,y)}$. I claim that the associated forcing $P_I$ makes
the ground model reals null. It will be enough to find a Polish measure space $Y,\mu$ and a Borel
set $B\subset X\times Y$ such that the vertical sections of $B$ are $\mu$-null and the horizontal
sections of the complement are $\gs$-porous. 

Consider the measure space for adding a random subset $U$ of the tree $T$ in such
a way that $t\in U:t\in T$ are mutually independent events and $\mu(t\in U)=1-p_{|t|}$. Let $U$ be the
generic random set. Use the definitions to show
that whenever $x\in X$ is a ground model point then only finitely many initial segments of the sequence $x$
are in the set $U$ (this happens because $\Pi_n p_n\neq 0$) but there are infinitely many initial segments
of $x$ which have an immediate successor in the set $U$ (this happens because $\Pi_n p_n^{k_n}=0$).
Now in the model $V[U]$ enumerate the set $U$ as $U=\{t_i:i\in\gw\}$ and for every number $j\in\gw$
let $A_j=\{x\in X\cap V:\forall i>j\ t_i\not\subset x\}$. The previous observations imply that
this set has porosity $\geq 1/2$ at each of its points and therefore is porous, and moreover
each point in $X\cap V$ is in one of the sets $A_j:j\in\gw$. So $V[U]\models
X\cap V$ is $\gs$-porous. The argument is completed by translating this conclusion to the existence
of the suitable Borel subset of the product $X\times Y$.
\end{example}

\section{Determinacy}

The purpose of this section is to give two similar determinacy arguments that were postponed in the previous sections.
They are of interest independently of the purposes of the present paper.

Suppose that $\mu$ is a pavement submeasure on some Polish space $X$.
derived form a countable set $U$ of Borel pavers, with weight function
$w$. Suppose $A\subset X$ is a set, and $\eps>0$ is a number. Consider a game $G_\eps(A)$ between players Adam
and Eve. In the game, Adam gradually builds a set $W\subset U$ of pavers with $w(W)\leq\eps$ and Eve
builds a point $x\in X$. Eve wins if $x\in A\setminus\bigcup W$ and I will want to prove

\begin{lemma}
\label{submeasurelemma}
If $\mu(A)<\eps$ then Adam has a winning strategy in the game $G_\eps(A)$ which in turn implies $\mu(A)\leq\eps$.
\end{lemma}

First I must specify the schedule for both players in the game. At round $n$, Adam must play a finite set $W_n\subset U$
such that $w(W)\leq\eps$ and if $n\in m$ then $W_n\subset W_m$ and $w(W_m)-w(W_n)\leq 2^{-n}$. The set $W$ is recovered
in the end as $W=\bigcup_nW_n$.
For Eve, fix a Borel bijection $f:\cantor\to X$. Eve will play bits $b_n\in 2$ and she is allowed to tread
water, that is, wait an arbitrary number of rounds before placing another nontrivial move. The point $x\in X$ is then
recovered as $x=f(\langle b_n:n\in\gw\rangle)$.

To prove the lemma, note that if $\mu(A)<\eps$ then Adam has a winning strategy in which he can ignore Eve's moves
altogether: he can produce a set $W$ of weight $<\eps$ so that $A\subset\bigcup W$, winning no matter what Eve plays.
On the other hand, suppose that $\mu(A)>\eps$ and $\gs$ is a strategy for Adam. I must produce a counterplay
against the strategy which results in a point $x\in A\setminus W$, that is, it ends with Eve's victory.

First, a bit of notation. If $\tau$ is a finite play observing the strategy $\gs$ let $W_\tau$ be the 
finite collection of sets
Adam put into his set $W\subset U$ so far, and let $V_\tau$ be the collection of sets the strategy
$\gs$ will dictate him to put into $W$ in the infinite play extending $\tau$ in which Eve makes
only trivial moves past $\tau$. Thus $W_\tau\subset V_\tau$ and $w(V_\tau\setminus W_\tau)\leq 2^{-|\tau|}$.
For a number $n\in\gw$ and a bit $b\in 2$ let also $\tau nb$ be the extension of the play $\tau$ in which
Eve makes only trivial moves except at the round $n$, which is also the final round of the play $\tau nb$,
at which she plays the bit $b$. Finally, for every play $\tau$ and a bit $b\in 2$ let $B_{\tau b}=\bigcap_n\bigcup
(V_{\tau nb}\setminus W_{\tau nb}$; clearly $\mu(B_{\tau b})=0$.

Now since the set $A\subset X$ is of submeasure $>\eps$, there is a point $x\in A\setminus (\bigcup V_0\cup\bigcup
\{B_{\tau b}:\tau$ is a finite play and $b\in 2\})$, because the latter set has $\mu$ submeasure at most $\eps$. 
It is now easy to inductively build a sequence $0=\tau_0\subset\tau_1\subset\dots$ of partial counterplays against the strategy
$\gs$ so that Eve builds the point $x$ and $x\notin\bigcup V_{\tau_n}$ for every number $n$.
Such a play clearly leads to Eve's victory.

There are now several corollaries.

\begin{corollary}
\label{submeasurecorollary1}
(ZF+AD) Every set has an analytic subset of the same submeasure. (ZFC+LC) Every universally Baire set has
an analytic subset of the same submeasure.
\end{corollary}

\begin{proof}
The second assertion follows from the first since suitable large cardinal assumptions imply that for every
universally Baire set $A\subset X$, the model $\lr[A]$ satisfies the Axiom of Determinacy. The first assertion
is an immediate corollary of Lemma~\ref{submeasurelemma}. Suppose AD holds, and $A\subset X$ is a set
of submeasure $\eps>0$. For every number $n\in\gw$, Adam does not have a winning strategy in the game $G_{\eps-2^{-n}}
(A)$, and therefore Eve must have a winning strategy $\gs_n$. Let $A_n\subset X$ be the set of all points $x\in X$
resulting from some Adam's counterplay against the strategy $\gs_n$. Then $A_n\subset A$ since the strategy
$\gs_n$ is winning, $A_n$ is an analytic set by its definition, and $\mu(A_n)\geq\eps-2^{-n}$ since the strategy
$\gs_n$ remains winning for Eve in the game $G_{\eps-2^{-n}}(A_n)$. Clearly, the set $\bigcup_nA_n\subset A$ is the required
analytic set.
\end{proof}

\begin{corollary}
\label{submeasurecorollary2}
(ZF+DC+AD+) The submeasure is continuous in increasing wellordered unions of uncountable cofinality.
\end{corollary}

\begin{proof}
Work in the theory ZF+DC+AD+.
Let $\kappa\in\theta$ be an uncountable regular cardinal. By a theorem of Steel \cite{jackson:square}, there is a set
$B\subset\R$ and a prewellordering $\prec$ on the set $B$ of length $\kappa$ such that every analytic
subset of the set $B$ meets only $<\kappa$ many classes of $\prec$. Suppose that $\langle C_\ga:\ga\in\kappa\rangle$
is an increasing union of sets and $\eps=\sup_\ga \mu(C_\ga)$. I must argue that $\mu(\bigcup_\ga C_\ga)=\eps$.
Consider the submeasure $\mu^*$ on $X\times\R$ given by $\mu^*(E)=\mu($projection of $E$ into the $X$ coordinate$)$.
It is immediate that this is a pavement submeasure derived from a countable
collection of Borel pavers. Consider the set $E\subset X\times\R$
given by $\langle x,r\rangle\in E\liff x\in C_\ga$ where $r\in B$ is a real number in the $\ga$-th class of the
prewellordering $\prec$. Clearly, $\mu^*(E)=\mu(\bigcup_\ga C_\ga)$. By the previous Corollary applied to
the submeasure $\mu^*$, there is an analytic subset $F\subset E$ of the same $\mu^*$ submeasure. The projection of the
set $F$ into the $\R$ coordinate meets less than $\kappa$ many classes, bounded by some ordinal $\gb\in\kappa$.
Then the projection of the set $F$ into the $X$ coordinate is a subset of the set $C_\gb$ and it has the
same $\mu$ submeasure as the set $\bigcup_\ga C_\ga$. Thus $\mu(C_\gb)=\mu(\bigcup_\ga C_\ga)$ as desired.
\end{proof}

\begin{corollary}
\label{submeasurecorollary3}
(ZF+DC+AD) Every set has a Borel subset of the same submeasure. (ZFC+LC) Every universally Baire set
has a Borel subset of the same submeasure.
\end{corollary}

\begin{proof}
There are two ways to argue here. The first is to show that in ZFC, every analytic set has a Borel subset of the same
submeasure and use Corollary~\ref{submeasurecorollary1}. This argument is interesting in its own right. 
Consider forcing with the partial order $P$ of analytic subsets of the Polish space $X$ with positive $\mu$ submeasure, ordered by
inclusion. Since the ideal $I$ of $\mu$-null sets is generated by Borel sets, an argument similar to
\cite{z:four} 4.17 shows that in the $P$-generic extension there is a real $\dot x_{gen}$
which belongs to all sets in the generic filter, falls out of analytic sets which are not in the generic filter,
and falls out of all $I$-small sets.
Let $A\in P$ be
an analytic $\mu$-positive set, and let $M$ be a countable elementary submodel of some large structure.
The proof of Theorem 7.4 of \cite{z:four} shows that the set $B=\{x\in X:\exists g\subset P\cap M$
an $M$-generic filter such that $A\in g\land x=\dot x_{gen}/g\}$ has $\mu$-submeasure equal to the set $A$.
Now $B\subset A$ since the condition $A$ forces in $P$ that $\dot x_{gen}\in\dot A$ and by forcing theorem then,
for every point $x\in B$, $M[x]\models x\in A$ and by absoluteness $x\in A$. Also the set $B$ is Borel:
it is in
one-to-one Borel correspondence with the $G_\gd$ set of all $M$-generic filters on $P\cap M$ containing
the condition $A\in P$. The Corollary follows. In retrospect of course the forcing $P_I$ is dense in the poset $P$.

The second way is to argue that in ZF, every analytic set is an increasing union
of $\gw_1$ many Borel sets, and with the assumption of AD the argument for Corollary~\ref{submeasurecorollary2}
shows that one of the Borel sets must have the same submeasure as the analytic set. 
\end{proof}

The treatment of strongly subadditive capacities is similar. Suppose that $c$ is a strongly subadditive capacity
on some Polish space $X$. Suppose $A\subset X$ is a set and $\eps>0$ is a real number. The setup of the game $G_\eps(A)$
is literally the same, with Adam playing basic open sets $W_n\subset X$ of capacity $\leq\eps$ coming from some
fixed countable basis for the space $X$ closed under finite unions. In particular, if $n\in m$ then 
$W_n\subset W_m$ and $c(W_m)-c(W_n)\leq
2^{-3n}$. The only miniscule change is in the last exponent, it is necessary
for purely arithmetical reasons. In the end, set $W=\bigcup_n W_n$. Again, the key point is

\begin{lemma}
$c(A)<\eps$ implies that Adam has a winning strategy in the game $G_\eps(A)$, which in turn implies
that $c(A)\leq\eps$.
\end{lemma}

As before, if $c(A)<\eps$ then $A\subset O$ for some open set $O\subset X$ of capacity $<\eps$, and Adam can easily
win by producing the set $W=O$, ignoring Eve's moves altogether. On the other hand, assume that $\gs$ is Adam's
strategy and $c(A)>\eps$; I must produce a counterplay against the strategy $\gs$ winning for Eve. To do this,
for an arbitrary number $n\in\gw$ argue that the union $V(n)$ of all open sets the strategy $\gs$ can produce
against Eve's counterplays with the first $n$ moves trivial and all other moves
nontrivial, has capacity $<\eps+2^{-n}$. Once this is done, Eve will just pick a number $n\in\gw$
such that $c(A)>\eps+2^{-n}$, a binary sequence $r\in 2^\gw$ such that $f(r)\in A\setminus V_n$, and she will win by first
waiting for $n$ moves and then producing the sequence $r$ without further hesitation.

The fact that $c(V_n)<\eps+2^{-n}$ immediately follows from an easy and useful claim.

\begin{claim}
Suppose that $c$ is a strongly subadditive capacity, $n\in\gw$, and 
$f:2^{<\gw}\to\cP(X)$ is a map such that $t\subset s$ implies $f(t)\subset f(s)$ and $c(f(s))-c(f(t))<2^{-2|t|-1-n}$.
Then $c(\bigcup_tf_t)\leq c(f(0))+2^{-n}$.
\end{claim}

\begin{proof}
For every number $k\in\gw$ let $Z_k=\bigcup_{t\in 2^k}f(t)$.

By induction on $k\in\gw$ prove that $c(Z_k)\leq c(f(0))+\gS_{l\in k}2^{-l-n}$. The case $k=0$ is clear.
Suppose this is known for some $k$. Then $Z_{k+1}=Z_k\cup\bigcup_{t\in 2^{k+1}}f(s)$ and for each sequence $s\in 2^{k+1}$ 
there is a subset $Y$
of $Z_k$ (namely $f(s\restriction k)$) such that $c(f(s))-c(Y)<2^{-2k-1-n}$. The strong subadditivity
of the capacity $c$
then implies that $c(Z_{k+1})<c(Z_k)+2^{k+1}2^{-2k-1-n}=c(Z_k)+2^{-k-n}$ as desired.

The conclusion of the claim then follows by the continuity of the capacity under increasing unions.
\end{proof}

The corollaries are similar to the case of the pavement submeasures.

\begin{corollary}
\label{capacitycorollary1}
(ZF+AD) Every set has an analytic subset of the same capacity. (ZFC+LC) Every universally Baire set has
an analytic subset of the same capacity.
\end{corollary}

\begin{corollary}
(ZF+DC+AD+) The capacity is continuous in increasing wellordered unions.
\end{corollary}

This is the same as in Corollary~\ref{submeasurecorollary1}
noting that the capacity is by definition continuous in increasing unions of
countable cofinality.

\begin{corollary}
\label{capacitycorollary3}
(ZF+AD) Every set has an $F_\gs$ subset of the same capacity. (ZFC+LC) Every universally Baire set
has an $F_\gs$ subset of the same capacity.
\end{corollary}

\begin{proof}
This immediately follows from Choquet's capacitability theorem and Corollary~\ref{capacitycorollary1}. 
\end{proof}

Comparing this with the proof
of Corollary~\ref{submeasurecorollary3}, it is interesting to note that I do not know whether 
in general the factor forcing $P_I$ is necessarily proper, where
$I$ is the ideal of sets of zero capacity. This nevertheless turned out to be true for every strongly subadditive capacity
for which I was able to verify it.

\bibliographystyle{plain}
\bibliography{shelah,odkazy,zapletal}

\end{document}